# Markov chain comparison[*]

**Martin Dyer**

*School of Computing
University of Leeds
Leeds LS2 9JT, UK
e-mail:* `dyer@comp.leeds.ac.uk`

**Leslie Ann Goldberg**

*Department of Computer Science
University of Warwick
Coventry CV4 7AL, UK
e-mail:* `leslie.goldberg@dcs.warwick.ac.uk`

**Mark Jerrum**

*School of Informatics
University of Edinburgh
Edinburgh EH9 3JZ, UK
e-mail:* `mrj@inf.ed.ac.uk`

**Russell Martin**

*Department of Computer Science
University of Liverpool
Liverpool L69 7ZF, UK
e-mail:* `R.Martin@csc.liv.ac.uk`

**Abstract:** This is an expository paper, focussing on the following scenario. We have two Markov chains, $\mathcal{M}$ and $\mathcal{M}'$. By some means, we have obtained a bound on the mixing time of $\mathcal{M}'$. We wish to compare $\mathcal{M}$ with $\mathcal{M}'$ in order to derive a corresponding bound on the mixing time of $\mathcal{M}$. We investigate the application of the comparison method of Diaconis and Saloff-Coste to this scenario, giving a number of theorems which characterize the applicability of the method. We focus particularly on the case in which the chains are not reversible. The purpose of the paper is to provide a catalogue of theorems which can be easily applied to bound mixing times.

## 1. Introduction

This expository paper focusses on *Markov chain comparison,* which is an important tool for determining the *mixing time* of a Markov chain.

[*] Partially supported by the EPSRC grant Discontinuous Behaviour in the Complexity of Randomized Algorithms

We are interested in the following scenario. We have two Markov chains, $\mathcal{M}$ and $\mathcal{M}'$. By some means, we have obtained a bound on the mixing time of $\mathcal{M}'$. We wish to compare $\mathcal{M}$ with $\mathcal{M}'$ in order to derive a corresponding bound on the mixing time of $\mathcal{M}$.

The foundation for the method lies in the comparison inequalities of Diaconis and Saloff-Coste [5, 6]. These inequalities give bounds on the eigenvalues of a reversible Markov chain in terms of the eigenvalues of a second chain. Similar inequalities were used by Quastel [17] in his study of the simple exclusion process on coloured particles.

The inequalities of Diaconis and Saloff-Coste provide the foundation for obtaining mixing-time bounds via comparison because there is a known close relationship between the mixing time of an ergodic reversible Markov chains and the eigenvalue which is second largest (in absolute value). This relationship was made explicit by Diaconis and Stroock [8] and Sinclair [19, proposition 1]. The latter is a discrete-time version of a proposition of Aldous [1].

Randall and Tetali [18] were the first to combine Diaconis and Saloff-Coste's inequalities with the inequalities relating eigenvalues to mixing times to derive a relationship between the mixing times of two chains. Their result [18, Proposition 4] applies to two ergodic reversible chains $\mathcal{M}$ and $\mathcal{M}'$ provided the eigenvalues satisfy certain restrictions (see the remarks following Theorem 10 below).

While the inequalities of Diaconis and Saloff-Coste are stated for reversible Markov chains, their proof does not use reversibility.[1] The Dirichlet forms correspond more closely to mixing times in the time-reversible case, but there is still some correspondence even without reversibility, as has been observed by Mihail [15] and Fill [10].

The primary purpose of our article is to pin down the applicability of the comparison method for non-reversible chains. This is done in Section 4. The main result (Theorem 24) is rather weaker than the corresponding theorem for reversible chains (Theorem 8) but we give examples (Observation 21 and the remark following Theorem 22) pointing out that the additional constraints are necessary.

Section 3 describes the comparison theorem for reversible chains. The main result (Theorem 8) is proved using exactly the method outlined by Randall and Tetali [18]. We feel that it is useful to provide a general theorem (Theorem 8) which applies to all reversible chains, including those that do not satisfy constraints on the eigenvalues. Diaconis and Saloff-Coste's method is sufficient for this task, provided the construction of the comparison is based on "odd flows" rather than just on flows. Observation 12 shows that the restriction that the flows be odd is necessary. The statement of Theorem 8 is deliberately general in terms of the parameters $\varepsilon$ and $\delta$, which are deliberately different to each other (unlikely the corresponding theorem in [18]). The reason for the generality is that the freedom to choose $\delta$ can lead to stronger results, as illustrated by

---

[1]For non-reversible Markov chains, the eigenvalues of the transition matrix are not necessarily real, but it is still possible to make sense of "spectral gap" as we shall see in Section 2.

Example 9. We have included a proof of Theorem 5, which is essentially Proposition 1(ii) of Sinclair [19] because, as far as we know, no proof appears in the literature. The theorem gives a lower bound on the mixing time of an ergodic reversible Markov chain in terms of its eigenvalues. A continuous-time version has been proved by Aldous [1].

Lemma 27 in Section 5 formalizes a technique that we have found useful in the past. In order to use the comparison inequalities of Diaconis and Saloff-Coste, one must construct a flow in which the congestion on an edge of the transition graph of the Markov chain is small. Lemma 27 shows that it is sometimes sufficient to construct a flow in which the congestion on a state is small.

Finally, we note that Section 5 of Randall and Tetali's paper [18] surveys other comparison methods which are not based on Diaconis and Saloff-Coste's inequalities. We will not repeat this survey, but refer the reader to [18].

## 2. Comparing Dirichlet forms

The following variation of Diaconis and Saloff-Coste's comparison method comes from [5, Section C]. It adapts an idea of Sinclair [19].

### *2.1. Definitions*

Let $\mathcal{M}$ be an ergodic (connected and aperiodic) Markov chain with transition matrix $P$, stationary distribution $\pi$, and state space $\Omega$. In this paper, the state space $\Omega$ will always be discrete and finite. We will assume that all Markov chains are discrete-time chains except where we indicate otherwise. Let $E(\mathcal{M})$ be the set of pairs of distinct states $(x, y)$ with $P(x, y) > 0$. Let $E^*(\mathcal{M})$ be the set of all pairs $(x, y)$ (distinct or not) with $P(x, y) > 0$. We will sometimes refer to the members of $E^*(\mathcal{M})$ as "edges" because they are the edges of the transition graph of $\mathcal{M}$. Define the *optimal Poincaré constant* of $\mathcal{M}$ by

$$\lambda_1(\mathcal{M}) = \inf_{\varphi:\Omega\to\mathbb{R}} \frac{\mathcal{E}_\mathcal{M}(\varphi,\varphi)}{\operatorname{var}_\pi \varphi},$$

where the infimum is over all non-constant functions $\varphi$ from $\Omega$ to $\mathbb{R}$ and the *Dirichlet form* is given by

$$\mathcal{E}_\mathcal{M}(\varphi,\varphi) = \frac{1}{2}\sum_{x,y\in\Omega} \pi(x)P(x,y)(\varphi(x)-\varphi(y))^2$$

and

$$\operatorname{var}_\pi \varphi = \sum_{x\in\Omega}\pi(x)(\varphi(x)-E_\pi\varphi)^2 = \frac{1}{2}\sum_{x,y\in\Omega}\pi(x)\pi(y)(\varphi(x)-\varphi(y))^2.$$

Let

$$\mathcal{F}_\mathcal{M}(\varphi,\varphi) = \frac{1}{2}\sum_{x,y\in\Omega}\pi(x)P(x,y)(\varphi(x)+\varphi(y))^2.$$

If $N$ is the size of $\Omega$ then let

$$\lambda_{N-1}(\mathcal{M}) = \inf_{\varphi:\Omega\to\mathbb{R}} \frac{\mathcal{F}_\mathcal{M}(\varphi,\varphi)}{\mathrm{var}_\pi\, \varphi},$$

where again the infimum is over all non-constant functions. When $\mathcal{M}$ is time reversible, the eigenvalues $1 = \beta_0 \geq \beta_1 \geq \cdots \geq \beta_{N-1}$ of the transition matrix $P$ are real. Then (see Facts 3 and 4 below) $\lambda_1(\mathcal{M})$ may be interpreted as the gap between $\beta_1$ and $1$, while $\lambda_{N-1}(\mathcal{M})$ is the gap between $\beta_{N-1}$ and $-1$. Although this explains the notation, the definitions of $\lambda_1$ and $\lambda_{N-1}$ make sense even for non-reversible Markov chains.

Suppose that $\mathcal{M}$ is an ergodic Markov chain on state space $\Omega$ with transition matrix $P$ and stationary distribution $\pi$, and that $\mathcal{M}'$ is another ergodic Markov chain on the same state space with transition matrix $P'$ and stationary distribution $\pi'$.

For every edge $(x,y) \in E^*(\mathcal{M}')$, let $\mathcal{P}_{x,y}$ be the set of paths from $x$ to $y$ using transitions of $\mathcal{M}$. More formally, let $\mathcal{P}_{x,y}$ be the set of paths $\gamma = (x = x_0, x_1, \ldots, x_k = y)$ such that

1. each $(x_i, x_{i+1})$ is in $E^*(\mathcal{M})$, and
2. each $(z,w) \in E^*(\mathcal{M})$ appears at most twice on $\gamma$.[2]

We write $|\gamma|$ to denote the length of path $\gamma$. So, for example, if $\gamma = (x_0, \ldots, x_k)$ we have $|\gamma| = k$. Let $\mathcal{P} = \cup_{(x,y)\in E^*(\mathcal{M}')} \mathcal{P}_{x,y}$.

An $(\mathcal{M}, \mathcal{M}')$-*flow* is a function $f$ from $\mathcal{P}$ to the interval $[0,1]$ such that for every $(x,y) \in E^*(\mathcal{M}')$,

$$\sum_{\gamma \in \mathcal{P}_{x,y}} f(\gamma) = \pi'(x) P'(x,y). \tag{1}$$

The flow is said to be an *odd* $(\mathcal{M}, \mathcal{M}')$-*flow* if it is supported by odd-length paths. That is, for every $\gamma \in \mathcal{P}$, either $f(\gamma) = 0$ or $|\gamma|$ is odd.

Let $r((z,w), \gamma)$ be the number of times that the edge $(z,w)$ appears on path $\gamma$. For every $(z,w) \in E^*(\mathcal{M})$, the *congestion* of edge $(z,w)$ in the flow $f$ is the quantity

$$A_{z,w}(f) = \frac{1}{\pi(z) P(z,w)} \sum_{\gamma \in \mathcal{P}:(z,w)\in\gamma} r((z,w),\gamma)\, |\gamma|\, f(\gamma).$$

The *congestion* of the flow is the quantity

$$A(f) = \max_{(z,w)\in E^*(\mathcal{M})} A_{z,w}(f).$$

### 2.2. Theorems

The following theorems are due to Diaconis and Saloff-Coste [5]. Theorem 1 is Theorem 2.3 of [5].

---

[2] This requirement is there for technical reasons — we want to ensure that $\mathcal{P}_{x,y}$ is finite.

**Theorem 1.** *Suppose that $\mathcal{M}$ is an ergodic Markov chain on state space $\Omega$ and that $\mathcal{M}'$ is another ergodic Markov chain on the same state space. If $f$ is an $(\mathcal{M}, \mathcal{M}')$-flow then for every $\varphi : \Omega \to \mathbb{R}$, $\mathcal{E}_{\mathcal{M}'}(\varphi, \varphi) \leq A(f)\mathcal{E}_{\mathcal{M}}(\varphi, \varphi)$.*

**Remark.** *The statement of Theorem 2.3 in [5] requires $\mathcal{M}$ and $\mathcal{M}'$ to be reversible, but the proof does not use this fact. In [5], $\mathcal{P}_{x,y}$ is defined to be the set of simple paths from $x$ to $y$. This is not an important restriction, because a flow $f$ can always be transformed into a flow $f'$ which is supported by simple paths and satisfies $A(f') \leq A(f)$. We prefer to use our definition so that an odd flow is a special case of a flow.*

**Theorem 2.** *Suppose that $\mathcal{M}$ is an ergodic Markov chain on state space $\Omega$ and that $\mathcal{M}'$ is another ergodic Markov chain on the same state space. If $f$ is an odd $(\mathcal{M}, \mathcal{M}')$-flow then for every $\varphi : \Omega \to \mathbb{R}$, $\mathcal{F}_{\mathcal{M}'}(\varphi, \varphi) \leq A(f)\mathcal{F}_{\mathcal{M}}(\varphi, \varphi)$.*

**Remark.** *Theorem 2.2 of [5] corresponds to the special case in which $\mathcal{P}_{x,y}$ contains a particular path $\gamma$ with $f(\gamma) = \pi'(x)P'(x, y)$. Again, the statement of the theorem requires the chains to be reversible, but the proof does not use this. The authors point out that their theorem can be generalized to the flow setting as we have done in Theorem 2. The same proof works.*

**Remark.** *Note that the definition of an $(\mathcal{M}, \mathcal{M}')$-flow requires us to route $\pi'(x)P'(x, y)$ units of flow between $x$ and $y$ for every pair $(x, y) \in E^*(\mathcal{M}')$. If we do not care whether the resulting flow is an odd flow we can dispense with the case $x = y$. For this case, we can consider the length-0 path $\gamma$ from $x$ to itself, and we can assign $f(\gamma) = \pi'(x)P'(x, x)$. The quantity $f(\gamma)$ contributes nothing towards the congestion of the flow (since it does not use any edges of $E^*(\mathcal{M})$).*

## 3. Comparing reversible Markov chains

### 3.1. Definitions

The *variation distance* between distributions $\theta_1$ and $\theta_2$ on $\Omega$ is

$$||\theta_1 - \theta_2|| = \frac{1}{2}\sum_i |\theta_1(i) - \theta_2(i)| = \max_{A \subseteq \Omega} |\theta_1(A) - \theta_2(A)|.$$

For an ergodic Markov chain $\mathcal{M}$ with transition matrix $P$ and stationary distribution $\pi$, and a state $x$, the mixing time from $x$ is

$$\tau_x(\mathcal{M}, \varepsilon) = \min\left\{t > 0 : ||P^{t'}(x, \cdot) - \pi(\cdot)|| \leq \varepsilon \text{ for all } t' \geq t\right\}.$$

In fact, $||P^t(x, \cdot) - \pi(\cdot)||$ is non-increasing in $t$, so an equivalent definition is

$$\tau_x(\mathcal{M}, \varepsilon) = \min\left\{t > 0 : ||P^t(x, \cdot) - \pi(\cdot)|| \leq \varepsilon\right\}.$$

Let

$$\tau(\mathcal{M}, \varepsilon) = \max_x \tau_x(\mathcal{M}, \varepsilon).$$

Let $\mathcal{M}$ be an ergodic Markov chain with transition matrix $P$, stationary distribution $\pi$, and state space $\Omega$. Let $N = |\Omega|$. Suppose that $P$ is reversible with respect to $\pi$. That is, every $x, y \in \Omega$ satisfies $\pi(x)P(x,y) = \pi(y)P(y,x)$. Then the eigenvalues of $P$ are real numbers and the maximum eigenvalue, 1, has multiplicity one. The eigenvalues of $P$ will be denoted as $\beta_0 = 1 > \beta_1 \geq \cdots \geq \beta_{N-1} > -1$. Let $\beta_{\max}(\mathcal{M}) = \max(\beta_1, |\beta_{N-1}|)$. The spectral representation of the transition matrix $P$ plays an important role in the results that follow. Let $D$ denote the diagonal matrix $D = \text{diag}(\pi_0^{1/2}, \ldots, \pi_{N-1}^{1/2})$. Then, using the definition of reversibility, it is easy to see that $A = DPD^{-1}$ is a symmetric matrix. Therefore, standard results from linear algebra (for example, see [11]) tell us there exists an orthonormal basis $\{e^{(i)'} : 0 \leq i \leq N-1\}$ of left eigenvectors of $A$, where $e^{(i)'}$ is an eigenvector corresponding to the eigenvalue $\beta_i$, i.e., $e^{(i)'}A = \beta_i e^{(i)'}$. We also have that $e_j^{(0)} = \pi_j^{1/2}$ for $j \in \{0, \ldots, N-1\}$. The important result we require (see [13, Section 3.2] or [20, Proof of Prop. 2.1] for a derivation) is that for $n \in \mathbb{N}$,

$$P^n(j,k) = \pi(k) + \sqrt{\frac{\pi_k}{\pi_j}} \sum_{i=1}^{N-1} \beta_i^n e_j^{(i)} e_k^{(i)} \qquad (2)$$

where the $P^n(j,k)$ are the $n$-step transition probabilities.

The following facts are well-known from linear algebra and follow from the "minimax" (or "variational") characterization of the eigenvalues (see [11], in particular the Rayleigh-Ritz and Courant-Fischer theorems).

**Fact 3.** *Let $\beta_0 = 1 > \beta_1 \geq \cdots \geq \beta_{N-1} > -1$ be the eigenvalues of the transition matrix of a reversible Markov chain $\mathcal{M}$. Then $1 - \beta_1 = \lambda_1(\mathcal{M})$.*

**Fact 4.** *Let $\beta_0 = 1 > \beta_1 \geq \cdots \geq \beta_{N-1} > -1$ be the eigenvalues of the transition matrix of a reversible Markov chain $\mathcal{M}$. Then $1 + \beta_{N-1} = \lambda_{N-1}(\mathcal{M})$.*

### 3.2. Lower bounds on mixing time

The following theorem is the same as Proposition 1(ii) of [19] by Sinclair (apart from a factor of 2). Sinclair's proposition is stated without proof. Aldous [1] proves a continuous-time version of Theorem 5. As far as we are aware, there is no published proof of the lower bound in discrete time so, for completeness, we provide one here based on Aldous's idea.

**Theorem 5.** *Suppose $\mathcal{M}$ is an ergodic reversible Markov chain. Then, for $\varepsilon > 0$,*

$$\tau(\mathcal{M}, \varepsilon) \geq \frac{\beta_{\max}(\mathcal{M})}{1 - \beta_{\max}(\mathcal{M})} \ln \frac{1}{2\varepsilon}.$$

*Proof of Theorem 5.* Let $P$ be the transition matrix of $\mathcal{M}$ and write the eigenvalues of $P$ as $\beta_0 = 1 > \beta_1 \geq \cdots \geq \beta_{N-1} > -1$. Let $\pi$ be the stationary distribution of $\mathcal{M}$. Let $A$ be the matrix defined in the spectral representation of $P$ in Section 3.1. Let $d(n) = \max_{j \in \Omega} ||P^n(j, \cdot) - \pi(\cdot)||$. We first give a lower bound on $d(2n)$.

Let $e^{(\max)}$ denote an eigenvector (of $A$) corresponding to $\beta_{\max}$. Since $e^{(\max)}$ is an eigenvector (and hence not identically zero), there exists some coordinate $j_0$ with $e^{(\max)}_{j_0} \neq 0$. Then, using (2) we find

$$
\begin{aligned}
d(2n) &= \max_{j \in \Omega} \|P^{2n}(j, \cdot) - \pi(\cdot)\| \\
&= \max_{j \in \Omega} \frac{1}{2} \sum_{k \in \Omega} \left| \sqrt{\frac{\pi(k)}{\pi(j)}} \sum_{i=1}^{N-1} \beta_i^{2n} e_j^{(i)} e_k^{(i)} \right| \\
&\geq \frac{1}{2} \sum_{k \in \Omega} \left| \sqrt{\frac{\pi(k)}{\pi(j_0)}} \sum_{i=1}^{N-1} \beta_i^{2n} e_{j_0}^{(i)} e_k^{(i)} \right| \\
&\geq \frac{1}{2} \left| \sum_{i=1}^{N-1} \beta_i^{2n} \left( e_{j_0}^{(i)} \right)^2 \right| = \frac{1}{2} \sum_{i=1}^{N-1} \beta_i^{2n} \left( e_{j_0}^{(i)} \right)^2 \\
&\geq \frac{1}{2} \beta_{\max}^{2n} \left( e_{j_0}^{(\max)} \right)^2.
\end{aligned}
$$

Using this lower bound, we have

$$\liminf_n \frac{d(2n)}{\beta_{\max}^{2n}} \geq \liminf_n \frac{\frac{1}{2} \beta_{\max}^{2n} \left( e_{j_0}^{(\max)} \right)^2}{\beta_{\max}^{2n}} = \frac{1}{2} \left( e_{j_0}^{(\max)} \right)^2 > 0. \tag{3}$$

Fix $\delta > 0$, and let $\tau^* = \tau(\mathcal{M}, \frac{\delta}{2})$. So, by definition $d(\tau^*) \leq \frac{\delta}{2}$.

For $s, t \in \mathbb{N}$ it is known that $d(s+t) \leq 2d(s)d(t)$. (For a proof of this fact see [2, Chapter 2].) Using this inequality we see that $d(2k\tau^*) \leq 2^{2k}(d(\tau^*))^{2k} \leq \delta^{2k}$. Then

$$
\begin{aligned}
\liminf_k \frac{d(2k\tau^*)}{\beta_{\max}^{2k\tau^*}} &\leq \liminf_k \delta^{2k} \beta_{\max}^{-2k\tau^*} \\
&= \liminf_k \exp\left(2k(\ln \delta - \tau^* \ln \beta_{\max})\right). \tag{4}
\end{aligned}
$$

If $\ln \delta - \tau^* \ln \beta_{\max} < 0$, the $\liminf$ in (4) is 0. This contradicts the lower bound in (3) (this lower bound also applies to the subsequence $d(2k\tau^*)/\beta_{\max}^{2k\tau^*}$ of $d(2n)/\beta_{\max}^{2n}$). Thus we conclude that $\ln \delta - \tau^* \ln \beta_{\max} \geq 0$, or $\tau^* \geq \frac{\ln(1/\delta)}{\ln(1/\beta_{\max})}$.

Finally, assuming $\beta_{\max} > 0$ (otherwise the theorem holds trivially),

$$\ln \frac{1}{\beta_{\max}} = \int_{\beta_{\max}}^1 \frac{dx}{x} \leq \int_{\beta_{\max}}^1 \frac{dx}{x^2} = \left[-\frac{1}{x}\right]_{\beta_{\max}}^1 = \frac{1 - \beta_{\max}}{\beta_{\max}}.$$

Combining this inequality with the previous one we obtain

$$\tau^* = \tau\left(\mathcal{M}, \frac{\delta}{2}\right) \geq \frac{\beta_{\max}}{1 - \beta_{\max}} \ln \frac{1}{\delta}.$$

Taking $\delta = 2\varepsilon$ gives the theorem. □

**Corollary 6.** *Suppose $\mathcal{M}$ is an ergodic reversible Markov chain. Then*

$$\tau\left(\mathcal{M}, \frac{1}{2e}\right) \geq \frac{\beta_{\max}(\mathcal{M})}{1 - \beta_{\max}(\mathcal{M})}.$$

### *3.3. Upper bounds on mixing time*

The following theorem is due to Diaconis and Stroock [8, Proposition 3] and to Sinclair [19, Proposition 1(i)].

**Theorem 7.** *Suppose that $\mathcal{M}$ is an ergodic reversible Markov chain with stationary distribution $\pi$. Then*

$$\tau_x(\mathcal{M}, \varepsilon) \leq \frac{1}{1 - \beta_{\max}(\mathcal{M})} \ln \frac{1}{\varepsilon \pi(x)}.$$

### *3.4. Combining lower and upper bounds*

In the following theorem we combine Diaconis and Saloff Coste's comparison method (Theorems 1 and 2) with upper bounds on mixing time (Theorem 7) and lower bounds on mixing time (Theorem 5) to obtain a comparison theorem for mixing times (Theorem 8). This combination was first provided by Randall and Tetali in Proposition 4 of [18]. We use the same reasoning as Randall and Tetali, though we consider odd flows in order to avoid assuming that the eigenvalues are non-negative.

**Theorem 8.** *Suppose that $\mathcal{M}$ is a reversible ergodic Markov chain with stationary distribution $\pi$ and that $\mathcal{M}'$ is another reversible ergodic Markov chain with the same stationary distribution. Suppose that $f$ is an odd $(\mathcal{M}, \mathcal{M}')$-flow. Then, for any $0 < \delta < \frac{1}{2}$,*

$$\tau_x(\mathcal{M}, \varepsilon) \leq A(f) \left[ \frac{\tau(\mathcal{M}', \delta)}{\ln(1/2\delta)} + 1 \right] \ln \frac{1}{\varepsilon \pi(x)}.$$

*In particular,*

$$\tau_x(\mathcal{M}, \varepsilon) \leq A(f) \left[ \tau\left(\mathcal{M}', \frac{1}{2e}\right) + 1 \right] \ln \frac{1}{\varepsilon \pi(x)}. \quad (5)$$

*Proof.* Let $N$ be the size of the state space.

$$\begin{aligned}
\tau_x(\mathcal{M}, \varepsilon) &\leq \frac{1}{1 - \beta_{\max}(\mathcal{M})} \ln \frac{1}{\varepsilon \pi(x)} && \text{(by Theorem 7)} \\
&= \max\left(\frac{1}{\lambda_1(\mathcal{M})}, \frac{1}{\lambda_{N-1}(\mathcal{M})}\right) \ln \frac{1}{\varepsilon \pi(x)} && \text{(by Facts 3 and 4)} \\
&\leq \max\left(\frac{A(f)}{\lambda_1(\mathcal{M}')}, \frac{A(f)}{\lambda_{N-1}(\mathcal{M}')}\right) \ln \frac{1}{\varepsilon \pi(x)} && \text{(by Theorems 1 and 2)} \\
&= A(f) \frac{1}{1 - \beta_{\max}(\mathcal{M}')} \ln \frac{1}{\varepsilon \pi(x)} && \text{(by Facts 3 and 4)} \\
&\leq A(f) \left[\frac{\tau(\mathcal{M}', \delta)}{\ln(1/2\delta)} + 1\right] \ln \frac{1}{\varepsilon \pi(x)} \\
&&& \text{(by Theorem 5, noting } \ln(1/2\delta) > 0\text{).}
\end{aligned}$$

□

**Remark.** *The proof of Theorem 8 goes through if $\mathcal{M}$ and $\mathcal{M}'$ have different stationary distributions as long as they have the same state space. However, an extra factor arises to account for the difference between $\operatorname{var}_\pi \varphi$ and $\operatorname{var}_{\pi'} \varphi$ in the application of Theorems 1 and 2. Not much is known about comparison of two chains with very different state spaces. However, there have been some successes. See [7].*

The freedom to choose $\delta$ in the statement of Theorem 8 is often useful, as we see in the following example, based on the "hard core gas" model.

**Example 9.** *Suppose $\mathcal{G}$ is some class of graphs of maximum degree $\Delta$. Let $G \in \mathcal{G}$ be an n-vertex graph, and let $\Omega_G$ denote the set of all independent sets in $G$. Let $\mathcal{M}'_G$ be an ergodic Markov chain on $\Omega_G$ with uniform stationary distribution. Of the transitions of $\mathcal{M}'_G$ we assume only that they are uniformly distributed and local. Specifically, to make a transition, a vertex $v$ of $G$ is selected uniformly at random; then the current independent set is randomly updated just on vertices within distance $r$ of $v$. We regard the radius $r$, as well as the degree bound $\Delta$, as a constant.*

*Now let $\mathcal{M}_G$ be another Markov chain that fits this description with $r = 0$. That is, it makes single-site updates according, say, to the heat-bath rule. Suppose we have proved that $\tau(\mathcal{M}'_G, \varepsilon) = O(n \log(n/\varepsilon))$. (This is a typical form of mixing time bound coming out of a coupling argument.) In this example, any reasonable choice of flow $f$ will have $A(f) = O(1)$: the canonical paths are of constant length, and a constant number of them flow through any transition of $\mathcal{M}$. (To know how the constant implicit in $O(1)$ depends on $\Delta$ and $r$, we'd need to have more details about the transitions of $\mathcal{M}'$, and be precise about the flow $f$.) Note that it is easy to arrange for $f$ to be an odd flow. Applying Theorem 8 with the default choice $\delta = 1/2e$ yields $\tau(\mathcal{M}_G, \varepsilon) = O(n^2 \log(n/\varepsilon))$, whereas setting $\delta$ optimally at $\delta = 1/n$, we achieve $\tau(\mathcal{M}_G, \varepsilon) = O(n^2 \log(1/\varepsilon))$, gaining a factor $\log n$.*

In the literature, the applications of Diaconis and Saloff-Coste's comparison method to mixing times are typically presented for the special case in which $\beta_{\max}(\mathcal{M})$ is the second-highest eigenvalue of the transition matrix of $\mathcal{M}$. In this case, it is not necessary for the flow $f$ to be an odd flow, so the proof of Theorem 8 gives the following.

**Theorem 10.** *Recall that $\beta_1$ is the second-highest eigenvalue of the transition matrix of $\mathcal{M}$. Theorem 8 holds with "odd $(\mathcal{M}, \mathcal{M}')$-flow", replaced by "$(\mathcal{M}, \mathcal{M}')$-flow", provided $\beta_{\max}(\mathcal{M}) = \beta_1$.*

**Remark.** *Theorem 10 is similar to Randall and Tetali's inequality [18, Proposition 4], which assumes that $\beta_{\max}(\mathcal{M})$ and $\beta_{\max}(\mathcal{M}')$ correspond to the second-highest eigenvalues of the relevant transition matrices, and that the latter is at least $1/2$.*

Since the restriction that $f$ be an odd flow is usually omitted from applications of comparison to mixing in the literature, it is worth considering the following example, which shows that the restriction is crucial for Theorem 8. The general idea underlying the example is simple. Let $\beta_0 = 1 > \beta_1 \geq \cdots \geq \beta_{N-1}$ be the eigenvalues of the transition matrix of $\mathcal{M}$. The eigenvalue $\beta_{N-1}$ is equal to $-1$ if $\mathcal{M}$ is periodic and is greater than $-1$ otherwise. If this eigenvalue is close to $-1$ then $\mathcal{M}$ is nearly periodic, and this slows down the mixing of $\mathcal{M}$. Let $\mathcal{M}'$ be the uniform random walk on the state space of $\mathcal{M}$. Clearly $\mathcal{M}'$ mixes in a single step, but we can construct a $(\mathcal{M}, \mathcal{M}')$-flow with low congestion as long as we take care to send flow along paths whose lengths are consistent with the (near) periodicity of $\mathcal{M}$.

**Example 11.** *Let $\Omega = \{a, b\}$. For a parameter $\delta \in (0, 1)$, let $P(a, b) = P(b, a) = 1 - \delta$ and $P(a, a) = P(b, b) = \delta$. The stationary distribution $\pi$ of $\mathcal{M}$ is uniform. Let $t$ be even. Then*

$$\|P^t(a, \cdot) - \pi(\cdot)\| \geq \Pr(X_t = a \mid X_0 = a) - \frac{1}{2} \geq (1 - \delta)^t - \frac{1}{2} \geq \frac{1}{2} - \delta t,$$

*so $\tau_a(\mathcal{M}, 1/4) \geq \lfloor 1/(2\delta) \rfloor$. Let $\mathcal{M}'$ be the uniform random walk on $\Omega$. The chain $\mathcal{M}'$ has stationary distribution $\pi$ and mixes in a single step.*

*Let $P$ be the transition matrix of $\mathcal{M}$ and let $P'$ be the transition matrix of $\mathcal{M}'$. We will now construct a $(\mathcal{M}, \mathcal{M}')$-flow $f$.*

*For the edge $(a, b) \in E^*(\mathcal{M}')$ let $\gamma$ be the length-1 path $(a, b)$ and assign $f(\gamma) = \pi(a)P'(a, b)$. Similarly, for the edge $(b, a)$, let $\gamma$ be the length-1 path $(b, a)$ and assign $f(\gamma) = \pi(b)P'(b, a)$. For the edge $(a, a) \in E^*(\mathcal{M}')$, let $\gamma$ be the path $a, b, a$ and assign $f(\gamma) = \pi(a)P'(a, a)$. Finally, for the edge $(b, b)$ let $\gamma$ be the path $b, a, b$ and assign $f(\gamma) = \pi(b)P'(b, b)$.*

*Note that $A_{a,a}(f) = A_{b,b}(f) = 0$. Also,*

$$A_{a,b}(f) = \frac{1}{\pi(a)P(a,b)} \left( \pi(a)P'(a,b) + 2\pi(a)P'(a,a) \right) = \frac{1}{1 - \delta} \left( \frac{1}{2} + 1 \right) = \frac{3}{2(1 - \delta)}.$$

*Similarly, $A_{b,a}(f) = 3/(2(1 - \delta))$. If we take $\delta \leq 1/2$, we have $A(f) \leq 3$. Then*

$$\tau_a(\mathcal{M}, \tfrac{1}{4}) \geq \left\lfloor \frac{1}{2\delta} \right\rfloor \frac{A(f)}{3} \frac{\tau(\mathcal{M}', \frac{1}{2e}) + 1}{2} \frac{\ln\left(1/\frac{1}{4}\pi(a)\right)}{\ln 8},$$

*and by making $\delta$ small, we can get as far away as we want from the inequality (5) in Theorem 8.*

Example 11 prompts the following observation.

**Observation 12.** *In general, for reversible ergodic Markov chains $\mathcal{M}$ and $\mathcal{M}'$ with the same state space and stationary distribution, the ratio between $\tau_x(\mathcal{M}, \varepsilon)$ and the quantity $\left[ \tau(\mathcal{M}', \frac{1}{2e}) + 1 \right] \ln \frac{1}{\varepsilon \pi(x)}$ from the right-hand-side of (5) can not be upper-bounded in terms of the congestion of an $(\mathcal{M}, \mathcal{M}')$-flow. (We know from Theorem 8 that such a bound is possible if we restrict attention to odd flows.)*

It is well-known (see, for example, Sinclair [19]) that the eigenvalues of the transition matrix $P$ of a Markov chain $\mathcal{M}$ are all non-negative if every state has a self-loop probability which is at least $1/2$. That is, the eigenvalues are non-negative if every state $x$ satisfies $P(x,x) \geq 1/2$. Thus, Theorem 10 applies to any such Markov chain $\mathcal{M}$. Observation 13 below shows that even weaker lower bounds on self-loop probabilities can routinely be translated into mixing-time inequalities without consideration of odd flows.

**Observation 13.** *Suppose that $\mathcal{M}$ is a reversible ergodic Markov chain with transition matrix $P$ and stationary distribution $\pi$ and that $\mathcal{M}'$ is another reversible ergodic Markov chain with the same stationary distribution. Suppose that $f$ is a $(\mathcal{M},\mathcal{M}')$-flow. Let $c = \min_x P(x,x)$, and assume $c > 0$. Then, for any $0 < \delta < 1/2$,*

$$\tau_x(\mathcal{M},\varepsilon) \leq \max\left\{ A(f)\left[\frac{\tau(\mathcal{M}',\delta)}{\ln(1/2\delta)} + 1\right], \frac{1}{2c} \right\} \ln \frac{1}{\varepsilon\pi(x)}.$$

*Proof.* Write $P = cI + (1-c)\widehat{P}$. Since the matrix $\widehat{P}$ is stochastic, its eigenvalues $\hat{\beta}_i$ all satisfy $|\hat{\beta}_i| \leq 1$. The relationship between the eigenvalues of $P$ and those of $\widehat{P}$ is simply $\beta_i = c + (1-c)\hat{\beta}_i$. In particular $\beta_{N-1} = c + (1-c)\hat{\beta}_{N-1} \geq -1 + 2c$.

By Theorem 5,

$$\tau(\mathcal{M}',\delta) \geq \frac{\beta_{\max}(\mathcal{M}')}{1 - \beta_{\max}(\mathcal{M}')} \ln \frac{1}{2\delta}$$

or, equivalently,

$$\beta_{\max}(\mathcal{M}') \leq \frac{\frac{\tau(\mathcal{M}',\delta)}{\ln(1/2\delta)}}{\frac{\tau(\mathcal{M}',\delta)}{\ln(1/2\delta)} + 1}.$$

By Fact 3,

$$\lambda_1(\mathcal{M}') \geq \frac{1}{\frac{\tau(\mathcal{M}',\delta)}{\ln(1/2\delta)} + 1}.$$

and hence, by Theorem 1,

$$\lambda_1(\mathcal{M}) \geq \frac{1}{A(f)\left(\frac{\tau(\mathcal{M}',\delta)}{\ln(1/2\delta)} + 1\right)}.$$

On the other hand, we know by Fact 4 and the lower bound on $\beta_{N-1}$ calculated above that $\lambda_{N-1}(\mathcal{M}) \geq 2c$. Thus

$$1 - \beta_{\max}(\mathcal{M}) = \min\{\lambda_1(\mathcal{M}), \lambda_{N-1}(\mathcal{M})\} \geq \min\left\{ \frac{1}{A(f)\left(\frac{\tau(\mathcal{M}',\delta)}{\ln(1/2\delta)} + 1\right)}, 2c \right\}.$$

The result follows immediately from Theorem 7. □

Suppose that $\mathcal{M}$ is a reversible Markov chain with transition matrix $P$. Let $\mathcal{M}_{\mathrm{ZZ}}$ be the Markov chain on the same state space with transition matrix $P_{\mathrm{ZZ}} = \frac{1}{2}(I + P)$. $\mathcal{M}_{\mathrm{ZZ}}$ is often referred to as the "lazy" version of $\mathcal{M}$. In the literature, it is common to avoid considering negative eigenvalues by studying the lazy chain $\mathcal{M}_{\mathrm{ZZ}}$ rather than the chain $\mathcal{M}$, using an inequality like the following.

**Observation 14.** *Suppose that $\mathcal{M}$ is a reversible ergodic Markov chain with stationary distribution $\pi$ and that $\mathcal{M}'$ is another reversible ergodic Markov chain with the same stationary distribution. Suppose that $f$ is a $(\mathcal{M}, \mathcal{M}')$-flow. Then*

$$\tau_x(\mathcal{M}_{ZZ}, \varepsilon) \leq 2A(f)[\tau(\mathcal{M}', \frac{1}{2e}) + 1] \ln \frac{1}{\varepsilon \pi(x)}.$$

*Proof.* Since the eigenvalues of $\mathcal{M}_{ZZ}$ are non-negative, $1 - \beta_{\max}(\mathcal{M}_{ZZ}) = \lambda_1(\mathcal{M}_{ZZ})$, which is equal to $\frac{1}{2}\lambda_1(\mathcal{M})$. The proof is now the same as that of Theorem 8. □

The approach of Observation 14 is somewhat unsatisfactory because it gives no insight about the mixing time of $\mathcal{M}$ itself. For example, we can applying Observation 14 to the chains $\mathcal{M}$ and $\mathcal{M}'$ from Example 11. This gives

$$\tau(\mathcal{M}_{ZZ}, 1/4) \leq 2 \frac{3}{2(1-\delta)}[1+1]\ln(8) = \frac{6}{1-\delta}\ln 8,$$

whereas we know from Example 11 that $\tau_a(\mathcal{M}, 1/4) \geq \lfloor 1/(2\delta) \rfloor$. Making $\delta \leq 1/4$ small, we can make $\mathcal{M}$ mix arbitrarily slowly, while $\mathcal{M}_{ZZ}$ mixes in at most 17 steps.

## 4. Comparison without reversibility

### 4.1. Definitions

Given a discrete-time ergodic Markov chain $\mathcal{M}$ on state space $\Omega$ with transition matrix $P$ and stationary distribution $\pi$, the *continuization* $\widetilde{\mathcal{M}}$ is defined as follows ([2, Chapter 2]). Let $Q$ be the *transition rate matrix* defined by $Q = P - I$. Then the distribution at time $t$ is $v \cdot \exp(Qt)$ where $v$ is the row vector corresponding to the initial distribution. (For a concise treatment of matrix exponentials, refer to Norris [16, Section 2.10].) The mixing time $\tau_x(\widetilde{\mathcal{M}}, \varepsilon)$ is thus

$$\tau_x(\widetilde{\mathcal{M}}, \varepsilon) = \inf\{t > 0 : ||v_x \cdot \exp(Qt') - \pi|| \leq \varepsilon \text{ for all } t' \geq t\},$$

where $v_x$ is the unit vector with a 1 in state $x$ and 0 elsewhere. Denote by $\widetilde{P}^t = \exp(Qt)$ the matrix of transition probabilities over a time interval of length $t$. A standard fact is the following (see Norris [16, Thm 2.1.1]).

**Lemma 15.** $\frac{d}{dt}\widetilde{P}^t = Q\widetilde{P}^t = \widetilde{P}^t Q$.

The *conductance* of a set $S$ of states of $\mathcal{M}$ is given by

$$\Phi_S(\mathcal{M}) = \frac{\sum_{i \in S}\sum_{j \in \overline{S}} \pi(i)P(i,j) + \sum_{i \in \overline{S}}\sum_{j \in S} \pi(i)P(i,j)}{2\pi(S)\pi(\overline{S})},$$

and the conductance of $\mathcal{M}$ is $\Phi(\mathcal{M}) = \min_S \Phi_S(\mathcal{M})$, where the min is over all $S \subset \Omega$ with $0 < \pi(S) < 1$.

Suppose $S$ is a subset of $\Omega$ with $0 < \pi(S) < 1$. Let $\chi_S$ be the indicator function for membership in $S$. That is, $\chi_S(x) = 1$ if $x \in S$ and $\chi_S(x) = 0$ otherwise. Then since $\operatorname{var}_\pi \chi_S = \pi(S)\pi(\overline{S})$, we have

$$\Phi_S(\mathcal{M}) = \frac{\mathcal{E}_\mathcal{M}(\chi_S, \chi_S)}{\operatorname{var}_\pi \chi_S}.$$

Thus, $\Phi(\mathcal{M})$ is the same as $\lambda_1(\mathcal{M})$ except that in $\Phi(\mathcal{M})$ we minimize over non-constant functions $\varphi : \Omega \to \{0, 1\}$ rather than over functions from $\Omega$ to $\mathbb{R}$. Thus we have

**Observation 16.** $\lambda_1(\mathcal{M}) \leq \Phi(\mathcal{M})$.

### 4.2. Lower bounds on mixing time

The analogue of Corollary 6 for the non-reversible case will be obtained by combining Theorems 17, 18 and 19 below. Theorem 17 is from Dyer, Frieze and Jerrum [9], and Theorem 18 is a continuous-time version of the same result.

**Theorem 17.** *Suppose that $\mathcal{M}$ is an ergodic Markov chain. Then*

$$\Phi(\mathcal{M}) \geq \frac{\frac{1}{2} - \frac{1}{2e}}{\tau(\mathcal{M}, \frac{1}{2e})}.$$

*Proof.* It is immediate from the symmetry in the definition of $\Phi_S(\mathcal{M})$ that $\Phi_S(\mathcal{M}) = \Phi_{\overline{S}}(\mathcal{M})$. Therefore, we can restrict the minimization in the definition of $\Phi(\mathcal{M})$ to sets $S$ with $0 < \pi(S) \leq \frac{1}{2}$. That is,

$$\Phi(\mathcal{M}) = \min_{S \subset \Omega, 0 < \pi(S) \leq \frac{1}{2}} \Phi_S(\mathcal{M}).$$

Also, since

$$\sum_{i \in \overline{S}} \sum_{j \in S} \pi(i) P(i, j) + \sum_{i \in \overline{S}} \sum_{j \in \overline{S}} \pi(i) P(i, j) = \sum_{i \in \overline{S}} \pi(i) \sum_{j \in \Omega} P(i, j) = \pi(\overline{S})$$

and

$$\pi(\overline{S}) = \sum_{i \in \Omega} \pi(i) \sum_{j \in \overline{S}} P(i, j) = \sum_{i \in S} \sum_{j \in \overline{S}} \pi(i) P(i, j) + \sum_{i \in \overline{S}} \sum_{j \in \overline{S}} \pi(i) P(i, j),$$

the two terms in the numerator in the definition of $\Phi_S(\mathcal{M})$ are identical and the definition of $\Phi_S(\mathcal{M})$ can be rewritten as follows.

$$\Phi_S(\mathcal{M}) = \frac{\sum_{i \in S} \sum_{j \in \overline{S}} \pi(i) P(i, j)}{\pi(S)\pi(\overline{S})}.$$

Let $\Phi'(\mathcal{M})$ be the asymmetrical version of conductance from [19]. Namely,

$$\Phi'(\mathcal{M}) = \min_{S \subset \Omega, 0 < \pi(S) \leq \frac{1}{2}} \Phi'_S(\mathcal{M}),$$

where
$$\Phi'_S(\mathcal{M}) = \Phi_S(\mathcal{M})\pi(\overline{S}).$$

In the process of proving Claim 2.3 in [9], Dyer, Frieze and Jerrum prove

$$\Phi'(\mathcal{M}) \geq \frac{\frac{1}{2} - \frac{1}{2e}}{\tau(\mathcal{M}, \frac{1}{2e})}.$$

The theorem follows since

$$\Phi(\mathcal{M}) = \min_{S \subset \Omega, 0 < \pi(S) \leq \frac{1}{2}} \Phi_S(\mathcal{M}) \geq \min_{S \subset \Omega, 0 < \pi(S) \leq \frac{1}{2}} \Phi_S(\mathcal{M})\pi(\overline{S}) = \Phi'(\mathcal{M}).$$

□

A very similar bound holds in continuous time. There does not seem to be a published proof of this result, so we provide one here, modelled on the proof in discrete time [9].

**Theorem 18.** *Suppose that $\mathcal{M}$ is an ergodic Markov chain. Then*

$$\Phi(\mathcal{M}) \geq \frac{\frac{1}{2} - \frac{1}{2e}}{\tau(\widetilde{\mathcal{M}}, \frac{1}{2e})}.$$

*Proof.* Suppose $\varphi$ is an arbitrary function $\varphi : \Omega \to \mathbb{R}$. The notation $\widetilde{P}^t\varphi$ denotes the function $\Omega \to \mathbb{R}$ defined by $[\widetilde{P}^t\varphi](x) = \sum_{y \in \Omega} \widetilde{P}^t(x,y)\varphi(y)$, for all $x \in \Omega$. Define $\varphi'_t : \Omega \to \mathbb{R}$ by $\varphi'_t(x) = \frac{d}{dt}\widetilde{P}^t\varphi(x)$ for all $x$. By Lemma 15, $\varphi'_t(x) = [Q\widetilde{P}\varphi](x) = [\widetilde{P}Q\varphi](x)$, where $[Q\widetilde{P}\varphi](x)$ is the function $\Omega \to \mathbb{R}$ defined by $[Q\widetilde{P}\varphi](x) = \sum_{y \in \Omega}(Q\widetilde{P})(x,y)\varphi(y)$ and $[\widetilde{P}Q\varphi](x)$ is defined similarly. Define $\theta : \mathbb{R}^+ \to \mathbb{R}^+$ by $\theta(t) = \|\varphi'_t\|_{\pi,1} = \sum_{x \in \Omega} \pi(x)|\varphi'_t(x)|$. Now observe that

$$|\widetilde{P}^t\varphi(x) - \varphi(x)| = \left|\int_0^t \varphi'_s(x)\,ds\right| \leq \int_0^t |\varphi'_s(x)|\,ds.$$

Multiplying by $\pi(x)$ and summing over $x$ we obtain

$$\|\widetilde{P}^t\varphi - \varphi\|_{\pi,1} \leq \int_0^t \|\varphi'_s\|_{\pi,1}\,ds = \int_0^t \theta(s)\,ds. \quad (6)$$

As before, denote by $\chi_S$ the indicator function of the set $S \subset \Omega$. We will show: (a) $\theta(t) \leq \theta(0)$, for all $t > 0$, and (b) if $\varphi = \chi_S$ then $\theta(0) \leq 2(\mathrm{var}_\pi \varphi)\Phi_S(\mathcal{M})$. It follows from these facts and (6) that

$$\|\widetilde{P}^t\varphi - \varphi\|_{\pi,1} \leq 2(\mathrm{var}_\pi \varphi)\Phi(\mathcal{M})t = 2\pi(S)\pi(\overline{S})\Phi(\mathcal{M})t, \quad (7)$$

where $\varphi = \chi_S$ and $S \subset \Omega$ is a set that minimises $\Phi_S(\mathcal{M})$, i.e., one for which $\Phi_S(\mathcal{M}) = \Phi(\mathcal{M})$. Now when $t = \tau(\widetilde{\mathcal{M}}, 1/2e)$, we know that $|\widetilde{P}^t\varphi(x) - \pi(S)| \leq$ $1/2e$ for all $x$. (Otherwise there would be a state $x \in \Omega$ for which $|\widetilde{P}^t\varphi(x) - \pi(S)| > 1/2e$. But then we would have $|\widetilde{P}^t(x,S) - \pi(S)| = |\widetilde{P}^t\varphi(x) - \pi(S)| > 1/2e$, contrary to the choice of $t$.) Assume without loss of generality that $\pi(S) \leq \frac{1}{2}$. Now

$$\sum_{x \in \Omega} \pi(x)[\widetilde{P}^t\varphi](x) = \sum_{x \in \Omega} \pi(x)\varphi(x) = \sum_{x \in S} \pi(x)\varphi(x),$$

so subtracting $\sum_{x \in S} \pi(x)[\widetilde{P}^t\varphi](x)$ from both sides (and using the fact that $\varphi$ is the indicator variable for $S$), we get

$$\sum_{x \in \Omega - S} \pi(x)([\widetilde{P}^t\varphi](x) - \varphi(x)) = \sum_{x \in S} \pi(x)(\varphi(x) - [\widetilde{P}^t\varphi](x)).$$

Then

$$\begin{aligned}
\|\widetilde{P}^t\varphi - \varphi\|_{\pi,1} &\geq 2\sum_{x \in S} \pi(x)|P^t\varphi(x) - \varphi(x)| \\
&\geq 2\sum_{x \in S} \pi(x)(1 - \pi(S) - 1/2e) \\
&\geq 2\sum_{x \in S} \pi(x)\left(\frac{1}{2} - \frac{1}{2e}\right) \\
&\geq 2\left(\frac{1}{2} - \frac{1}{2e}\right)\pi(S)\pi(\overline{S}),
\end{aligned}$$

which combines with (7) to yield the claimed result. To complete the proof, we need to verify facts (a) and (b).

By Lemma 15,

$$\begin{aligned}
\varphi'_t(x) &= \frac{d}{dt}\widetilde{P}^t\varphi(x) = \sum_{y \in \Omega} Q(x,y)[\widetilde{P}^t\varphi](y) \\
&= \left(\sum_{y \in \Omega} P(x,y)[\widetilde{P}^t\varphi](y)\right) - [\widetilde{P}^t\varphi](x) \\
&= \sum_{y \in \Omega} P(x,y)\left([\widetilde{P}^t\varphi](y) - [\widetilde{P}^t\varphi](x)\right).
\end{aligned}$$

In particular, if $\varphi = \chi_S$ and $t = 0$,

$$\theta(0) = \|\varphi'_0\|_{\pi,1} \leq \sum_{x,y \in \Omega} \pi(x)P(x,y)|\varphi(y) - \varphi(x)| = 2\mathcal{E}_{\mathcal{M}}(\varphi,\varphi) = 2(\mathrm{var}_\pi \varphi)\Phi_S(\mathcal{M}),$$

by definition of $\Phi_S$, which is fact (b). Fact (a) follows from the following sequence

of (in)equalities:

$$\begin{aligned}
\theta(t) = \|\varphi'_t\|_{\pi,1} = \sum_x \pi(x) \left|\varphi'_t(x)\right| &= \sum_x \pi(x) |[\widetilde{P}^t Q\varphi](x)| \\
&= \sum_x \pi(x) \left|\sum_y \widetilde{P}^t(x,y)[Q\varphi](y)\right| \\
&\leq \sum_x \pi(x) \sum_y \widetilde{P}^t(x,y) \left|[Q\varphi](y)\right| \\
&= \sum_y |[Q\varphi](y)| \sum_x \pi(x)\widetilde{P}^t(x,y) \\
&= \sum_y \pi(y) \left|[Q\varphi](y)\right| \\
&= \theta(0).
\end{aligned}$$

$\square$

The following theorem is known as Cheeger's inequality.

**Theorem 19.** *Suppose that $\mathcal{M}$ is an ergodic Markov chain. Then $\lambda_1(\mathcal{M}) \geq \Phi(\mathcal{M})^2/8$.*

*Proof.* We will reduce to the reversible case, in which Cheeger's inequality is well-known. Let $P$ be the transition matrix of $\mathcal{M}$ and let $\pi$ be its stationary distribution. Let $\widehat{\mathcal{M}}$ be the Markov chain with transition matrix $\widehat{P}(x,y) = \frac{1}{2}\bigl(P(x,y) + \pi(y)/\pi(x)P(y,x)\bigr)$. Now for any $\varphi : \Omega \to \mathbb{R}$

$$\begin{aligned}
\mathcal{E}_{\widehat{\mathcal{M}}}(\varphi,\varphi) &= \frac{1}{2} \sum_{x,y} \pi(x)\widehat{P}(x,y)(\varphi(x)-\varphi(y))^2 \\
&= \frac{1}{4} \sum_{x,y} \pi(x) P(x,y)(\varphi(x)-\varphi(y))^2 + \frac{1}{4} \sum_{x,y} \pi(y) P(y,x)(\varphi(x)-\varphi(y))^2 \\
&= \mathcal{E}_{\mathcal{M}}(\varphi,\varphi).
\end{aligned}$$

This implies both $\lambda_1(\mathcal{M}) = \lambda_1(\widehat{\mathcal{M}})$ and $\Phi(\mathcal{M}) = \Phi(\widehat{\mathcal{M}})$ since these are just minimisations of $\mathcal{E}_{\widehat{\mathcal{M}}}(\varphi,\varphi)$ and $\mathcal{E}_{\mathcal{M}}(\varphi,\varphi)$ over $\varphi$ (recall the remark just before Observation 16).

Note that $\widehat{\mathcal{M}}$ is time-reversible since

$$\pi(x)\widehat{P}(x,y) = \frac{1}{2}\pi(x)P(x,y) + \frac{1}{2}\pi(y)P(y,x) = \pi(y)\widehat{P}(y,x).$$

Now let $\Phi'(\widehat{\mathcal{M}})$ be the asymmetrical conductance from the proof of Theorem 17. Since $\widehat{\mathcal{M}}$ is time-reversible, the eigenvalues of $\widehat{P}$ are real numbers $\beta_0 = 1 > \beta_1 \geq \cdots \geq \beta_{N-1} > -1$ and from Fact 3 we have $1 - \beta_1 = \lambda_1(\widehat{\mathcal{M}})$. Now by Lemma 2.4 of Sinclair [20] we have $\lambda_1(\mathcal{M}) = \lambda_1(\widehat{\mathcal{M}}) \geq \Phi'(\widehat{\mathcal{M}})^2/2$. Also,

$$\Phi(\widehat{\mathcal{M}}) = \min_{S \subset \Omega, 0 < \pi(S) \leq \frac{1}{2}} \Phi_S(\widehat{\mathcal{M}}) \leq \min_{S \subset \Omega, 0 < \pi(S) \leq \frac{1}{2}} \Phi_S(\widehat{\mathcal{M}}) 2\pi(\overline{S}) = 2\Phi'(\widehat{\mathcal{M}}),$$

so $\Phi'(\widehat{\mathcal{M}}) \geq \Phi(\widehat{\mathcal{M}})/2 = \Phi(\mathcal{M})/2$. $\square$

Combining Theorems 17, 18 and 19 we get the following analogue of Corollary 6.

**Corollary 20.** *Suppose that $\mathcal{M}$ is an ergodic Markov chain. Then*

$$\lambda_1(\mathcal{M}) \geq \frac{(\frac{1}{2} - \frac{1}{2e})^2}{8\,\tau(\mathcal{M}, \frac{1}{2e})^2}$$

*and*

$$\lambda_1(\mathcal{M}) \geq \frac{(\frac{1}{2} - \frac{1}{2e})^2}{8\,\tau(\widetilde{\mathcal{M}}, \frac{1}{2e})^2}.$$

Note that the first lower bound is a function of the discrete mixing time $\tau(\mathcal{M}, \frac{1}{2e})$ and the second is a function of the continuous mixing time $\tau(\widetilde{\mathcal{M}}, \frac{1}{2e})$. Otherwise, the bounds are the same.

Corollary 20 seems quite weak compared to Corollary 6 because of the squaring of the mixing time in the denominator. It turns out that our bound cannot be improved for the general (non-reversible) case.

**Observation 21.** *There is an ergodic Markov chain $\mathcal{M}$ with $\lambda_1(\mathcal{M}) \in O(1/\tau(\mathcal{M}, \frac{1}{2e})^2)$.*

*Proof.* Let $\mathcal{M}$ be the Markov chain described in Section 2 of [4]. This chain has state space
$$\Omega = \{-(n-1), -(n-2), \ldots, -1, 0, 1, \ldots, n\}.$$
The transition matrix $P$ is defined as follows. If $j \equiv i+1 \pmod{2n}$ then $P(i,j) = 1 - 1/n$. Also, if $j \equiv -i \pmod{2n}$ then $P(i,j) = 1/n$. The stationary distribution $\pi$ is uniform on the $2n$ states. Diaconis et al. [4, Theorem 1] show that $\tau(\mathcal{M}, \frac{1}{2e}) \in O(n)$. We show that $\lambda_1(\mathcal{M}) \in O(1/n^2)$, implying $\lambda_1(\mathcal{M}) \in O(1/\tau(\mathcal{M}, \frac{1}{2e})^2)$.

Let $\varphi$ be the function given by $\varphi(i) = |i|$. We show that $\mathcal{E}_\mathcal{M}(\varphi, \varphi) \in O(1)$ and $\operatorname{var}_\pi \varphi \in \Omega(n^2)$.

First, the transitions from $i$ to $-i$ preserve $\varphi$, so to calculate $\mathcal{E}_\mathcal{M}$, we need only consider the edges from $i$ to $i+1$ (over which $\varphi$ differs by 1).

$$\mathcal{E}_\mathcal{M}(\varphi, \varphi) = \frac{1}{2} \sum_{x \in \Omega} \frac{1}{2n}(1 - 1/2n).$$

To calculate $\operatorname{var}_\pi \varphi$, we observe that

$$\mathbb{E}_\pi \varphi = \frac{1}{2n}\left(\sum_{i=1}^{n-1} i + \sum_{i=1}^{n} i\right) = \frac{1}{2n}\left(2\frac{n(n-1)}{2} + n\right) = \frac{n}{2},$$

so $(\mathbb{E}_\pi \varphi)^2 = n^2/4$. Also

$$\mathbb{E}_\pi(\varphi^2) = \frac{1}{2n}\left(n^2 + 2\sum_{i=1}^{n-1} i^2\right) = \frac{1}{2n}\left(n^2 + 2\frac{(n-1)n(2(n-1)+1)}{6}\right) = \frac{2n^2+1}{6}.$$

So $\operatorname{var}_\pi \varphi = \mathbb{E}_\pi(\varphi^2) - (\mathbb{E}_\pi \varphi)^2 = \frac{n^2+2}{12}$. □

## 4.3. Upper bounds on mixing time

We now give a continuous-time analogue of Theorem 7 for the general (non-reversible) case. Our proof follows the treatment in Jerrum [12] on pages 63+ followed by 55-56, with a few details filled in.

**Theorem 22.** *Suppose that $\mathcal{M}$ is an ergodic Markov chain with stationary distribution $\pi$. Then*

$$\tau_x(\widetilde{\mathcal{M}}, \varepsilon) \leq \frac{1}{2\lambda_1(\mathcal{M})} \ln\left(\frac{1}{\varepsilon^2 \pi(x)}\right).$$

*Proof.* Let $\Omega$ be the state space of $\mathcal{M}$. Let $P$ be the transition matrix of $M$, let $Q = P - I$ be the transition rate of $\widetilde{\mathcal{M}}$ and $\widetilde{P}^t = \exp(Qt)$ as in Section 4.1. Let $\varphi$ be any function from $\Omega$ to $\mathbb{R}$ with $\mathbb{E}_\pi \varphi = 0$. If $x$ is chosen from a distribution then $[\widetilde{P}^t \varphi](x)$ is a random variable. Note that $\mathbb{E}_\pi[\widetilde{P}^t \varphi] = 0$.

Using Lemma 15, we have

$$\frac{d}{dt} \operatorname{var}_\pi[\widetilde{P}^t \varphi] = \frac{d}{dt} \sum_{x \in \Omega} \pi(x)\left([\widetilde{P}^t \varphi](x)\right)^2$$

$$= 2 \sum_{x \in \Omega} \pi(x)[\widetilde{P}^t \varphi](x) \frac{d}{dt}[\widetilde{P}^t \varphi](x)$$

$$= 2 \sum_{x \in \Omega} \pi(x)[\widetilde{P}^t \varphi](x)[Q\widetilde{P}^t \varphi](x)$$

$$= 2 \sum_{x,y \in \Omega} \pi(x)Q(x,y)[\widetilde{P}^t \varphi](x)[\widetilde{P}^t \varphi](y)$$

$$= 2 \sum_{x,y \in \Omega} \pi(x)P(x,y)[\widetilde{P}^t \varphi](x)[\widetilde{P}^t \varphi](y) - 2\sum_{x \in \Omega} \pi(x)[\widetilde{P}^t \varphi](x)[\widetilde{P}^t \varphi](x)$$

$$= -2\mathcal{E}_\mathcal{M}(\widetilde{P}^t \varphi, \widetilde{P}^t \varphi)$$

$$= -2 \frac{\mathcal{E}_\mathcal{M}(\widetilde{P}^t \varphi, \widetilde{P}^t \varphi)}{\operatorname{var}_\pi[\widetilde{P}^t \varphi]} \operatorname{var}_\pi[\widetilde{P}^t \varphi]$$

$$\leq -2\lambda_1(\mathcal{M}) \operatorname{var}_\pi[\widetilde{P}^t \varphi].$$

Now let $w$ denote $\operatorname{var}_\pi[\widetilde{P}^t \varphi]$ and consider the differential equation $\frac{d}{dt}w \leq -2\lambda_1(\mathcal{M})w$. Solving this we get $\frac{dw}{w} \leq -2\lambda_1(\mathcal{M})dt$ so $\ln w \leq -2\lambda_1(\mathcal{M})t + c$ and $w \leq \exp(c - 2\lambda_1(\mathcal{M})t)$. Plugging in $t = 0$ we get $w \leq (\operatorname{var}_\pi \varphi) \exp(-2\lambda_1(\mathcal{M})t)$. For a subset $A \subseteq \Omega$ define $\varphi : \Omega \to \mathbb{R}$ by

$$\varphi(x) = \begin{cases} 1 - \pi(A) & \text{if } x \in A, \\ -\pi(A) & \text{otherwise.} \end{cases}$$

Note that $\mathbb{E}_\pi \varphi = 0$ and $\operatorname{var}_\pi \varphi \leq 1$ so

$$\operatorname{var}_\pi[\widetilde{P}^t \varphi] \leq \exp(-2\lambda_1(\mathcal{M})t).$$

Set
$$t = \frac{1}{2\lambda_1(\mathcal{M})} \ln\left(\frac{1}{\varepsilon^2 \pi(x)}\right)$$

so $\operatorname{var}_\pi[\widetilde{P}^t \varphi] \leq \varepsilon^2 \pi(x)$. Now $\operatorname{var}_\pi[\widetilde{P}^t \varphi] \geq \pi(x)([\widetilde{P}^t \varphi](x))^2$ so

$$\varepsilon \geq |\widetilde{P}^t \varphi(x)| = \left|\sum_y \exp(Qt)(x,y)\varphi(y)\right|$$
$$= \left|\left(\sum_{y \in A} \exp(Qt)(x,y)\right) - \pi(A)\right|$$
$$= \left|\Pr(\widetilde{X}_t \in A) - \pi(A)\right|,$$

for any $A$. So $\tau_x(\widetilde{\mathcal{M}}, \varepsilon) \leq t$. □

**Remark.** *If $\mathcal{M}$ is not reversible, it may not always be possible to get a discrete-time version of Theorem 22. In fact, we cannot always upper-bound $\tau_x(\mathcal{M}, \varepsilon)$ even as a function of both $1/\lambda_1(\mathcal{M})$ and $1/\lambda_{N-1}(\mathcal{M})$. Here is an example. Let $\mathcal{M}$ be the Markov chain which deterministically follows a directed cycle of length 3. Both $\lambda_1(\mathcal{M})$ and $\lambda_{N-1}(\mathcal{M})$ are bounded above 0. To see this, let $\varphi$ be a function from $\Omega$ to $\mathbb{R}$ with variance 1. Note that $\mathcal{E}_\mathcal{M}(\varphi, \varphi)$ and $\mathcal{F}_\mathcal{M}(\varphi, \varphi)$ are non-zero. But $\mathcal{M}$ is not ergodic (so certainly we cannot upper-bound its mixing time!). Let $\mathcal{M}'$ denote the uniform random walk on $K_3$. There is a low-congestion odd $(\mathcal{M}, \mathcal{M}')$-flow $f$, so $1/\lambda_1(\mathcal{M}) \leq A(f)/\lambda_1(\mathcal{M}')$ and $1/\lambda_{N-1}(\mathcal{M}) \leq A(f)/\lambda_{N-1}(\mathcal{M}')$. These inequalities do not give upper bounds on $\tau_x(\mathcal{M}, \varepsilon)$ because, while they rule out length-2 periodicities in $\mathcal{M}$, they do not rule out higher periodicities.*

Let $R(\mathcal{M})$ be the time-reversal of $\mathcal{M}$ with transition matrix $R(P)$ given by

$$R(P)(x,y) = \frac{\pi(y)}{\pi(x)} P(y,x).$$

Consider the chain $R(\mathcal{M})\mathcal{M}$ which does one step of $R(\mathcal{M})$ followed by one step of $\mathcal{M}$ during each transition. Here is a discrete-time companion to Theorem 22. This is based on Theorem 2.1 of Fill [10]. This idea (bounding convergence in terms of the Dirichlet form of $R(\mathcal{M})\mathcal{M}$) is also in [14].

**Theorem 23.** *Suppose that $\mathcal{M}$ is an ergodic Markov chain with stationary distribution $\pi$. Then*

$$\tau_x(\mathcal{M}, \varepsilon) \leq \frac{1}{\lambda_1(R(\mathcal{M})\mathcal{M})} \ln\left(\frac{1}{\varepsilon^2 \pi(x)}\right).$$

*Proof.* Let $\varphi$ be a function from $\Omega \to \mathbb{R}$ with $\mathbb{E}_\pi \varphi = 0$. The following equality, due to Mihail [15], is Proposition 2.3 of [10].

$$\mathcal{E}_{R(\mathcal{M})\mathcal{M}}(\varphi,\varphi) = \frac{1}{2}\sum_{x,y,z}\pi(x)R(P)(x,y)P(y,z)(\varphi(x)-\varphi(z))^2$$
$$= \frac{1}{2}\sum_{x,y,z}\pi(y)P(y,x)P(y,z)(\varphi(x)-\varphi(z))^2$$
$$= \operatorname{var}_\pi \varphi - \sum_{x,y,z}\pi(y)P(y,x)P(y,z)\varphi(x)\varphi(z)$$
$$= \operatorname{var}_\pi \varphi - \sum_y \pi(y)\sum_x P(y,x)\varphi(x)\sum_z P(y,z)\varphi(z)$$
$$= \operatorname{var}_\pi \varphi - \operatorname{var}_\pi(P\varphi).$$

This gives

$$\operatorname{var}_\pi(P\varphi) = \operatorname{var}_\pi \varphi - \mathcal{E}_{R(\mathcal{M})\mathcal{M}}(\varphi,\varphi)$$
$$= \left(1 - \frac{\mathcal{E}_{R(\mathcal{M})\mathcal{M}}(\varphi,\varphi)}{\operatorname{var}_\pi \varphi}\right)\operatorname{var}_\pi \varphi$$
$$\leq \left(1 - \lambda_1(R(\mathcal{M})\mathcal{M})\right)\operatorname{var}_\pi \varphi$$

so

$$\operatorname{var}_\pi(P^t\varphi) \leq \left(1 - \lambda_1(R(\mathcal{M})\mathcal{M})\right)^t \operatorname{var}_\pi \varphi \leq \exp(-t\lambda_1(R(\mathcal{M})\mathcal{M}))\operatorname{var}_\pi \varphi.$$

Then we can finish as in the proof of Theorem 22. □

### 4.4. Combining lower and upper bounds

The following theorem follows immediately from Theorem 22, Theorem 1 and Corollary 20.

**Theorem 24.** *Suppose that $\mathcal{M}$ is an ergodic Markov chain with stationary distribution $\pi$ and that $\mathcal{M}'$ is another ergodic Markov chain with the same stationary distribution. Suppose that $f$ is an $(\mathcal{M},\mathcal{M}')$-flow. Then*

$$\tau_x(\widetilde{\mathcal{M}},\varepsilon) \leq 4A(f)\frac{\tau(\widetilde{\mathcal{M}'},\frac{1}{2e})^2}{(\frac{1}{2}-\frac{1}{2e})^2}\ln\left(\frac{1}{\varepsilon^2 \pi(x)}\right)$$

*and*

$$\tau_x(\widetilde{\mathcal{M}},\varepsilon) \leq 4A(f)\frac{\tau(\mathcal{M}',\frac{1}{2e})^2}{(\frac{1}{2}-\frac{1}{2e})^2}\ln\left(\frac{1}{\varepsilon^2 \pi(x)}\right).$$

As in Corollary 20, the first inequality gives an upper bound in terms of the continuous mixing time $\tau(\widetilde{\mathcal{M}'},\frac{1}{2e})$. The second inequality is the same except that the upper bound is in terms of the discrete mixing time $\tau(\mathcal{M}',\frac{1}{2e})$.

If we use Theorem 23 instead of Theorem 22, we get

**Theorem 25.** *Suppose that $\mathcal{M}$ is an ergodic Markov chain with stationary distribution $\pi$ and that $\mathcal{M}'$ is another ergodic Markov chain with the same stationary distribution. Suppose that $f$ is an $(R(\mathcal{M})\mathcal{M}, \mathcal{M}')$-flow. Then*

$$\tau_x(\mathcal{M}, \varepsilon) \le A(f) \frac{8\tau(\mathcal{M}', \frac{1}{2e})^2}{\left(\frac{1}{2} - \frac{1}{2e}\right)^2} \ln\left(\frac{1}{\varepsilon^2 \pi(x)}\right).$$

We can do better if $\mathcal{M}'$ is known to be reversible. Using Corollary 6, we get

**Theorem 26.** *Suppose that $\mathcal{M}$ is an ergodic Markov chain with stationary distribution $\pi$ and that $\mathcal{M}'$ is another reversible ergodic Markov chain with the same stationary distribution. Suppose that $f$ is an $(\mathcal{M}, \mathcal{M}')$-flow. Then*

$$\tau_x(\widetilde{\mathcal{M}}, \varepsilon) \le \frac{A(f)}{2} \left[\tau\left(\mathcal{M}', \frac{1}{2e}\right) + 1\right] \ln\left(\frac{1}{\varepsilon^2 \pi(x)}\right).$$

## 5. Comparison and state congestion

This Section generalizes an idea that we used in [3]. In order to use the comparison inequalities of Diaconis and Saloff-Coste, one must construct a flow in which the congestion on an edge of the transition graph of the Markov chain is small. The following lemma shows that it is sometimes sufficient to construct a flow in which the congestion on a state is small.

Suppose that $f$ is an $(\mathcal{M}, \mathcal{M}')$-flow. The congestion of state $z$ in $f$ is the quantity

$$B_z(f) = \frac{1}{\pi(z)} \sum_{\gamma \in \mathcal{P}: z \in \gamma} |\gamma| f(\gamma).$$

Let $B(f) = \max_{z \in \Omega} B_z(f)$. Let

$$\kappa(f) = \max_{(z,w): A_{z,w}(f) > 0} \left(\sum_{x \in \Omega} \min\left\{P(z, x), R(P)(w, x)\right\}\right)^{-1}.$$

**Lemma 27.** *Suppose that $f$ is an $(\mathcal{M}, \mathcal{M}')$-flow. Then there is an $(\mathcal{M}, \mathcal{M}')$-flow $f'$ with $A(f') \le 8\kappa(f)B(f)$.*

*Proof.* Without loss of generality, we will assume that the flow $f$ is supported by simple paths. That is $f(\gamma) = 0$ if the path $\gamma$ has a repeated vertex. See the remark after Theorem 1.

Let $p(z, w, x)$ denote $\min\{P(z, x), R(P)(w, x)\}$ and let $\delta(z, w) = \sum_{x \in \Omega} p(z, w, x)$. Construct $f'$ as follows. For every path $\gamma = (x_0, \ldots, x_k)$, route the $f(\gamma)$ units of flow from $x_0$ to $x_k$ along a collection of paths of length $2k$. In particular, spread the flow along $\gamma$ from $x_i$ to $x_{i+1}$ as follows. For each $x \in \Omega$, route $p(x_i, x_{i+1}, x) f(\gamma) / \delta(x_i, x_{i+1})$ of this flow along the route $x_i, x, x_{i+1}$.

First we check that $f'$ is an $(\mathcal{M}, \mathcal{M}')$-flow. Note that if $p(x_i, x_{i+1}, x) > 0$ then both $P(x_i, x) > 0$ and $R(P)(x_{i+1}, x) > 0$ so (since $\pi(x_{i+1})$ and $\pi(x)$ are assumed to be nonzero) $P(x, x_{i+1}) > 0$. We conclude that the edges used by $f'$ are edges of $E^*(\mathcal{M})$. Also, each edge appears at most twice, as required, since $f$ is simple.

Now we bound the congestion of $f'$. Let $(z,w)$ be an edge in $E^*(\mathcal{M})$. By definition, the congestion of edge $(z,w)$ in $f'$ is

$$A_{z,w}(f') = \frac{1}{\pi(z)P(z,w)} \sum_{\gamma' \in \mathcal{P}:(z,w)\in\gamma'} r((z,w),\gamma')\,|\gamma'|\,f'(\gamma')$$
$$\leq \frac{2}{\pi(z)P(z,w)} \sum_{\gamma' \in \mathcal{P}:(z,w)\in\gamma'} |\gamma'|\,f'(\gamma'). \tag{8}$$

But the flow $f'$ was constructed by "spreading" the flow $f(\gamma)$ on each $\gamma \in \mathcal{P}$ over a number of paths $\gamma'$ with $|\gamma'| = 2|\gamma|$ as described above. Thus, the right-hand-side of (8) is at most

$$\frac{2}{\pi(z)P(z,w)} \left( \sum_{y\in\Omega} \sum_{\gamma\in\mathcal{P}:(z,y)\in\gamma} 2|\gamma| \frac{f(\gamma)p(z,y,w)}{\delta(z,y)} \right.$$
$$\left. + \sum_{y\in\Omega} \sum_{\gamma\in\mathcal{P}:(y,w)\in\gamma} 2|\gamma| \frac{f(\gamma)p(y,w,z)}{\delta(y,w)} \right)$$
$$\leq \frac{4\kappa(f)}{\pi(z)P(z,w)} \left( \sum_{y\in\Omega} \sum_{\gamma\in\mathcal{P}:(z,y)\in\gamma} |\gamma|\,f(\gamma)p(z,y,w) \right.$$
$$\left. + \sum_{y\in\Omega} \sum_{\gamma\in\mathcal{P}:(y,w)\in\gamma} |\gamma|\,f(\gamma)p(y,w,z) \right)$$
$$\leq \frac{4\kappa(f)}{\pi(z)P(z,w)} \left( \sum_{y\in\Omega} \sum_{\gamma\in\mathcal{P}:(z,y)\in\gamma} |\gamma|\,f(\gamma)P(z,w) \right.$$
$$\left. + \sum_{y\in\Omega} \sum_{\gamma\in\mathcal{P}:(y,w)\in\gamma} |\gamma|\frac{f(\gamma)\pi(z)P(z,w)}{\pi(w)} \right)$$
$$\leq 4\kappa(f)\bigl(B_z(f) + B_w(f)\bigr).$$

□